\documentclass[12pt]{article}
\usepackage{amsmath}
\usepackage{amsfonts}
\usepackage{amssymb}
\usepackage{polski}

\input{tcilatex}
\begin{document}

\begin{center}
\textbf{Homomorphisms from Functional Equations: The Goldie Equation, II. }

\textbf{by}

\textbf{N. H. Bingham and A. J. Ostaszewski.}

\bigskip
\end{center}

\textbf{Abstract.} This first of three sequels to \textsl{Homomorphisms from
Functional equations: The Goldie equation} [Ost2] by the second author --
the second of the resulting quartet -- starts from the Goldie functional
equation arising in the general regular variation of our joint paper
[BinO5]. We extend the work there in two directions. First, we algebraicize
the theory, by systematic use of certain groups -- the \textit{Popa groups}
arising in earlier work by Popa, and their relatives the \textit{Javor groups%
}. Secondly, we extend from the original context on the real line to
multi-dimensional (or infinite-dimensional) settings.

\bigskip


\bigskip

\noindent \textbf{Classification}: 26A03, 26A12, 33B99, 39B22, 62G32.

\bigskip

\noindent \textbf{1. Introduction. }The \textit{Goldie functional equation} $%
(GFE)$ in its simplest form, involving as unknowns a primary function $K$
called a \textit{kernel} and an \textit{auxiliary} $g$, both \textit{%
continuous}, reads%
\begin{equation}
K(x+y)=K(x)+g(x)K(y).  \tag{$GFE$}
\end{equation}%
We encounter a more general version of $(GFE)$ below, a special case of a 
\textit{Levi-Civita equation}. The real-valued version above is closely
related to the better-known \textit{Go\l \k{a}b-Schinzel functional equation}%
\begin{equation}
\eta (x+y\eta (x))=\eta (x)\eta (y),  \tag{$GS$}
\end{equation}

It emerged most clearly in [BinO5] in the investigation of functions of
regular variation, where $(GFE)$ is key -- see \S 2 below, that that
equation is best studied by reference to \textit{Popa groups}. These involve
a group structure on $\mathbb{R}$ first introduced by Popa [Pop], defined by
the binary operation%
\begin{equation*}
x\circ y:=x+y\eta (x),
\end{equation*}%
which enables $(GS)$ to be restated as homomorphy of $\mathbb{G}_{\eta }^{+}(%
\mathbb{R}):=\{x:\eta (x)>0\}$ with the multiplicative group of positive
reals. Its generalization below to $\mathbb{R}^{d},$ has%
\begin{equation*}
\eta (x)\equiv 1+\rho (x)
\end{equation*}%
with $\rho (.)$ linear on $\mathbb{R}^{d}.$ With the induced Euclidean
topology, $\mathbb{G}_{\rho }(\mathbb{R}^{d})=\mathbb{G}_{1+\rho (.)}^{+}(%
\mathbb{R}^{d})$ is an open subspace of $\mathbb{R}^{d},$ so by the argument
in Hewitt and Ross [HewR, 15.18], for $\lambda _{d}$ Lebesgue measure, the
Popa Haar-measure on $\mathbb{G}_{\rho }(\mathbb{R}^{d})$ is (as in [BinO5])
proportional to%
\begin{equation*}
\frac{\lambda _{d}(\mathrm{d}x)}{1+\rho (x)}.
\end{equation*}%
This enables the identification of Fourier transforms, for instance for $%
\mathbb{G}_{\rho }(\mathbb{R})$ with $\rho \in (0,\infty ),$ 
\begin{equation*}
\hat{f}(\gamma )=\int_{\mathbb{G}_{\rho }}f(u)\gamma (u_{\rho }^{-1})(1+\rho
)\frac{du}{1+\rho u}\qquad (\gamma \in \mathbb{R}),
\end{equation*}%
where the characters take the form $u\mapsto e^{i\gamma \log (1+\rho u)}$
with $\gamma \in \mathbb{R}$ and $u_{\rho }^{-1}$ denotes inversion in the
group $\mathbb{G}_{\rho }(\mathbb{R})$.

It was noticed in [Ost2], again in the context of $\mathbb{R}$, that $(GFE)$
itself can be equivalently formulated as a homomorphy between a pair of Popa
groups on $\mathbb{R}$.

In this paper we develop radial properties of multivariate Popa groups in
order to characterize \textit{Popa homomorphisms} -- homomorphisms between
Popa groups.

Regular variation in one dimension (widely used in analysis, probability and
elsewhere -- cf. [Bin2]) explores the ramifications of limiting relations
such as%
\begin{equation}
f(\lambda x)/f(x)\rightarrow K(\lambda )\equiv \lambda ^{\gamma } 
\tag{$Kar_{\times }$}
\end{equation}%
or its additive variant, more thematic here:%
\begin{equation}
f(x+u)-f(x)\rightarrow K(u)\equiv \kappa u  \tag{$Kar_{+}$}
\end{equation}%
[BinGT, Ch. 1], and%
\begin{equation}
\lbrack f(x+u)-f(x)]/h(x)\rightarrow K(u)\equiv (u^{\gamma }-1)/\gamma 
\tag{$BKdH$}
\end{equation}%
(Bojani\'{c} \& Karamata, de Haan, [BinGT, Ch. 3]). Beurling regular
variation similarly explores the ramifications of relations such as%
\begin{equation}
\varphi (x+t\varphi (x))/\varphi (x)\rightarrow 1\text{ or }\eta (t) 
\tag{$Beu$}
\end{equation}%
[BinGT, \S\ 2.11] and [Ost1]. The underlying Popa structure lies disguised
in the limit function $\eta (t),$ which takes the form $1+\gamma t$ for $%
t>-1/\gamma .$

For background and applications, see the standard work [BinGT] and e.g.
[BinO1-6], [Bin1,2,3]. Both theory and applications prompt the need to work
in higher dimensions, finite or infinite. This is the ultimate motivation
for the present paper.

\bigskip

\textbf{2. The multivariate Goldie functional equation. }For $X$ a real
topological vector space, write $\langle u\rangle _{X}$ for the \textit{%
linear span} of $u\in X$ (to be differentiated from the use of $\langle
u\rangle _{\rho }$ below for $\rho $ in the dual of $X$). Following [Ost1]
call a function $\varphi :X\rightarrow \mathbb{R}$ \textit{self-equivarying}
over $X,$ $\varphi \in SE_{X},$ if for each $u\in X$ both $\varphi (tu)=O(t)$
and 
\begin{equation*}
\varphi (tu+v\varphi (tu))/\varphi (tu)\rightarrow \eta _{u}^{\varphi
}(v)\qquad (v\in \langle u\rangle _{X},t\rightarrow \infty )
\end{equation*}%
locally uniformly in $v.$ This appeals to the underlying uniformity
structure on $X$ generated by the neighbourhoods of the origin. As in [Ost1]
(by restriction to the linear span $\langle u\rangle _{X}$) the limit
function $\eta =\eta _{u}^{\varphi }$ satisfies $(GS)$ for $x,y\in \langle
u\rangle _{X}.$ When the limit function $\eta _{u}$ is continuous, one of
the forms it may take is%
\begin{equation*}
\eta _{u}(x)=1+\rho _{u}x\qquad (x\in \langle u\rangle _{X})
\end{equation*}%
for some $\rho _{u}\in \mathbb{R}$, the alternative form being $\eta
(x)=\max \{1+\rho _{u}x,0\}.$ A closer inspection of the proof in [Ost1]
shows that in fact the restriction $x,y\in \langle u\rangle _{X}$ placed on $%
(GS)$ above is unneccessary. Consequently, one may apply the Brillou\"{e}%
t-Dhombres-Brzd\k{e}k theorem [BriD, Prop. 3], [Brz1, Th. 4], on the
continuous solutions of $(GS)$ with $\eta :X\rightarrow \mathbb{R}$, to
infer that $\eta $ here takes the form%
\begin{equation*}
\eta (x)=1+\rho (x)\qquad (x\in X),
\end{equation*}%
for some continuous linear functional $\rho :X\rightarrow \mathbb{R}$, the
alternative form being $\eta (x)=\max \{1+\rho (x),0\}$. On this matter, see
also [Bar], [BriD], [Brz1]; cf. [Chu2,3], the former cited in detail below.
(For the same conclusion under assumptions such as radial continuity, or
Christensen measurability, see [Jab1,2], and [Brz2] under boundedness on a
non-meagre set.)

Below we study the implications of replacing $\rho _{u}$ in $\eta _{u}$ by a
continuous linear function $\rho (x).$ For this we now need to extend the
definition of \textit{general regular variation} [BinO5] from the real line
to a multivariate setting. For real topological vector spaces $X,Y,$ a
function $f:X\rightarrow Y$ is $\varphi $-\textit{regularly varying} for $%
\varphi \in SE_{X}$ relative to the (auxiliary)\textit{\ norming }function $%
h:X\rightarrow \mathbb{R}$ if the\textit{\ kernel} function $K$ below is
well defined for all $x\in X$ by%
\begin{equation}
K(x):=\lim_{t\rightarrow \infty }[f(tx+x\varphi (tx))-f(tx)]/h(tx)\qquad
(x\in X).  \tag{$GRV$}
\end{equation}%
For later use, we note the underlying \textit{radial dependence}: for $u\in
X $ put%
\begin{equation*}
K_{u}(x):=\lim_{s\rightarrow \infty }[f(su+x\varphi (su))-f(su)]/h(su)\qquad
(x\in \langle u\rangle _{X}).
\end{equation*}%
Writing $x=\xi u$ with $\xi >0$ and $s:=t\xi >0,$ 
\begin{eqnarray*}
K(x) &=&K(\xi u)=\lim_{t\rightarrow \infty }f(t\xi u+x\varphi (t\xi
u))-f(t\xi u)]/h(t\xi u) \\
&=&\lim_{s\rightarrow \infty }f(su+x\varphi (su))-f(su)]/h(su)=K_{u}(x).
\end{eqnarray*}%
So here $K_{u}=K|\langle u\rangle _{X},$ as $K(\xi u)=K_{u}(\xi u).$

\bigskip

We work radially: above with half-lines $(0,\infty )$ and below with those
of the form $(-1/\rho ,\infty )$ for $\rho >0$ (on $\langle u\rangle _{X}$
with context determining $u$) and $(-\infty ,\infty )$ when $\rho =0,$ see
[BinO5]. Proposition 1 below identifies the emergence of functional
equations satisfied by the kernel function $K:X\rightarrow Y$ and by its
other auxiliary $g$ defined below. The latter, once $\eta ^{\varphi }$ is
identified in the continuous context (for which see again [Ost1]), as above,
yields a multivariate form of $(GS).$ Given the natural association of the
auxilary to the Goldie equation, its defining multiplicative equation has
`dual citizenship', being both a special case of $GFE$ (take logarithms!)
and a partially pexiderized variant of $(GS),$ for which see [Chu1], [Jab1].

\bigskip

\noindent \textbf{Proposition 1.} \textit{Let }$h$\textit{\ and }$\varphi
\in SE_{X}$\textit{\ be such that the limit}%
\begin{equation*}
g(x):=\lim_{t\rightarrow \infty }h(tx+x\varphi (tx))/h(tx)\qquad (x\in X)
\end{equation*}%
\textit{exists. Then the kernel }$K:X\rightarrow Y$\textit{\ in }$(GRV)$ 
\textit{satisfies the Goldie functional equation:} 
\begin{equation}
K(x+\eta ^{\varphi }(x)y)=K(x)+g(x)K(y)  \tag{$GFE$}
\end{equation}%
\textit{for }$y\in \langle x\rangle _{X}.$ \textit{Furthermore, }$g$\textit{%
\ satisfies }$(GFE)$\textit{\ in the alternative form} 
\begin{equation}
g(x+\eta ^{\varphi }(x)y)=g(x)g(y)\qquad (y\in \langle x\rangle _{X}). 
\tag{$GS/GFE_{\times }$}
\end{equation}

\bigskip

\noindent \textbf{Proof. }Fix $x$ and $y.$ Writing $s=s_{x}:=t+\varphi (tx),$
so that $sx=tx+x\varphi (tx),$ 
\begin{equation*}
\frac{f(tx+(x+y)\varphi (tx))-f(tx)}{h(tx)}\hspace{3.5in}
\end{equation*}%
\begin{eqnarray*}
&=&\frac{f(sx+y[\varphi (tx)/\varphi (sx)]\varphi (sx))-f(sx)}{h(sx)}\cdot 
\frac{h(tx+x\varphi (tx))}{h(tx)} \\
&&+\frac{f(tx+x\varphi (tx))-f(tx)}{h(tx)}.
\end{eqnarray*}
Here $\varphi (sx)/\varphi (tx)=\varphi (tx+x\varphi (tx))/\varphi
(tx)\rightarrow \eta (x).$ Passage to the limit yields $(GFE),$ since $%
\varphi (tx)=O(t).$ The final assertion is similar but simpler. \hfill $%
\square $

\bigskip

We will achieve a characterization of the kernel function $K$ by identifying
the dependence between the different \textit{radial restrictions} $K|\langle
u\rangle _{X}.$

\bigskip

\textbf{3. Popa-Javor circle groups and their radial subgroups. }For a real
topological vector space $X$ and a continuous linear function $\rho
:X\rightarrow \mathbb{R},$ the associated function%
\begin{equation*}
\varphi (x)=\eta _{\rho }(x):=1+\rho (x)
\end{equation*}%
satisfies $(GS)$, as may be routinely checked. The associated circle
operation $\circ _{\rho }:$%
\begin{equation*}
x\circ _{\rho }y=x+y\varphi (x)=x+y+\rho (x)y
\end{equation*}%
(which gives for $\rho (x)=I(x)\equiv x$ and $X=\mathbb{R}$ the \textit{%
circle operation }of ring theory: cf. [Jac, II.3], [Coh, 3.1], and [Ost2, \S %
2.1] for the historical background) is due to Popa in 1965 on the line and
by Javor in 1968 in a vector space ([Pop], [Jav], cf. [BinO4]). It is
associative, as noted in [Jav]. As in [BinO5] we need the open sets%
\begin{equation*}
\mathbb{G}_{\rho }=\mathbb{G}_{\rho }(X):=\{x\in X:\eta _{\rho }(x)=1+\rho
(x)>0\}.
\end{equation*}%
Note that if $x,y\in \mathbb{G}_{\rho }$, then $x\circ _{\rho }y\in \mathbb{G%
}_{\rho }$, as 
\begin{equation*}
\eta _{\rho }(x\circ _{\rho }y)=\eta _{\rho }(x)\eta _{\rho }(y)>0.
\end{equation*}%
\textbf{Definition. }We refer to 
\begin{equation*}
\mathbb{G}_{\rho }^{\ast }=\mathbb{G}_{\rho }^{\ast }(X):=\{x\in X:\eta
_{\rho }(x)\neq 0\}
\end{equation*}%
as the \textit{Javor group }since, as Javor [Jav] shows, the set is a group
under $\circ _{\rho }.$ The Javor result remains true under the additional
restriction $\eta _{\rho }(y)>0,$ as we are about to verify in Theorem J\
below. Thus, likewise, we refer to 
\begin{equation*}
\mathbb{G}_{\rho }=\mathbb{G}_{\rho }(X):=\{x\in X:\eta _{\rho }(x)>0\}
\end{equation*}%
as a \textit{Popa group} under $\circ _{\rho }.$

\bigskip

\noindent \textbf{Theorem J }(after Javor [Jav])\textbf{.}\textit{\ For }$X$%
\textit{\ a topological vector space and }$\rho ~:~X~\rightarrow ~\mathbb{R}$%
\textit{\ a continuous linear function,} $(\mathbb{G}_{\rho }(X),\circ
_{\rho })$ \textit{is a group.}

\bigskip

\noindent \textbf{Proof. }This is routine, and one argues just as in [Jav],
but must additionally check preservation of the positivity of $\eta _{\rho }$
on $\mathbb{G}_{\rho }$. Here $0\in \mathbb{G}_{\rho }$ and is the neutral
element; the inverse of $x\in \mathbb{G}_{\rho }$ is $x_{\rho
}^{-1}:=-x/(1+\rho (x)),$ which is in $\mathbb{G}_{\rho }$ since $1=\eta
_{\rho }(0)=\eta _{\rho }(x)\eta _{\rho }(x_{\rho }^{-1}),$ so that $\eta
_{\rho }(x_{\rho }^{-1})>0.$ \hfill $\square $

\bigskip

\noindent \textbf{Definitions.} 1. For $u\in \mathbb{G}_{\rho }(X)$, put%
\begin{equation*}
\langle u\rangle _{\rho }:=\langle u\rangle _{X}\cap \mathbb{G}_{\rho
}(X)=\{tu:\eta _{\rho }(tu)=1+t\rho (u)>0,t\in \mathbb{R}\}.
\end{equation*}%
(If $\rho (u)\neq 0,$ then $\langle u\rangle _{\rho }=\{tu:t>-1/\rho (u)\},$
which is a half-line in $\langle u\rangle _{X}$; otherwise $\langle u\rangle
_{\rho }=\langle u\rangle _{X}.$ Note that $\mathbb{G}_{\rho }(X)$ is an
affine half-space in $X.$)

Given the context, the notation $\langle u\rangle _{\rho }$ will not clash
with that of $\langle u\rangle _{X}.$

\noindent 2. For $K$ with domain $\mathbb{G}_{\rho }(X)$ we will write $%
K_{u}=K|\langle u\rangle _{\rho }$. (This too will not clash with the radial
notation of \S 2.)

\bigskip

\noindent \textbf{Lemma.} \textit{The one-dimensional subgroup} $\langle
u\rangle _{\rho }$ \textit{is an abelian subgroup of }$\mathbb{G}_{\rho }(X)$%
\textit{\ isomorphic with }$\mathbb{G}_{\rho (u)}(\mathbb{R}).$

\bigskip

\noindent \textbf{Proof. }We check closure under multiplication and
inversion. For $s,t\in \mathbb{R},$ as before $\varphi (su\circ _{\rho
}tu)=\varphi (su)\varphi (tu)>0;$ also, writing $r(tu)$ for the $\rho $%
-inverse, $\varphi (r(tu))>0$ for $\varphi (tu)>0,$ as $1=\varphi
(0)=\varphi (tu\circ _{\rho }r(tu))=\varphi (tu)\varphi (r(tu)).$ Further,
since%
\begin{equation*}
su\circ _{\rho }tu=su+tu+st\rho (u)u=(s\circ _{\rho (u)}t)u,
\end{equation*}%
the operation $\circ _{\rho }$ is abelian on $\langle u\rangle _{\rho }.$%
\hfill $\square $

\bigskip

\noindent \textbf{Remark. }Despite the lemma above, unless $\rho \equiv 0$
or $X=\mathbb{R},$ the group $\mathbb{G}_{\rho }(X)$ itself is non-abelian.
(In the commutative case, except when $X=\mathbb{R},$ one may select $x\neq
0 $ with $\rho (x)=0;$ then $x\rho (y)=y\rho (x)=0$ and so $\rho (y)=0$ for
all $y.$) We return to this matter in detail in Theorem 2 below.

\bigskip

\noindent \textbf{Definition. }Say that a subgroup $H$ of $\mathbb{G}_{\rho
}(X)$ is \textit{radial }if $H\subseteq \langle u\rangle _{\rho }$ for some $%
u\in H$.

\bigskip

Theorem 1 below concerns radial subgroups. The assumption there on $\Sigma $
is effectively that all its radial subgroups are closed and dense in
themselves. Key to the proof is the observation that if $1+\rho (u)<0$, then
a fortiori $1+\rho (-u)=1-\rho (u)>0,$ i.e. if $u\notin \langle u\rangle
_{\rho },$ then its negative $-u\in \langle u\rangle _{\rho }$ and likewise
its $\mathbb{G}_{\rho }(X)$-inverse $(-u)_{\rho }^{-1}\in \langle u\rangle
_{\rho }$.

\bigskip

\noindent \textbf{Theorem 1.}\newline
\textit{Radial subgroups of Popa groups are Popa. That is, for }$\Sigma $%
\textit{\ a subgroup of} $\mathbb{G}_{\rho }(X)$ \textit{with }$\langle
u\rangle _{\rho }\subseteq \Sigma $\textit{\ for each }$u\in \Sigma :$%
\begin{equation*}
\Sigma =\mathbb{G}_{\rho }(\langle \Sigma \rangle _{X}).
\end{equation*}

\bigskip

\noindent \textbf{Proof. }With $\langle \Sigma \rangle $ the linear span, $%
\Sigma \subseteq \mathbb{G}_{\rho }(\langle \Sigma \rangle _{X})$ follows
from $\Sigma \subseteq \langle \Sigma \rangle _{X}$, as $\Sigma $ and $%
\mathbb{G}_{\rho }(\langle \Sigma \rangle _{X})$ are subgroups of $\mathbb{G}%
_{\rho }(X)$.

For the converse, we first show that $\alpha x+\beta y\in \Sigma $ for $%
x,y\in \Sigma $ and scalars $\alpha ,\beta $ whenever $\alpha x+\beta y\in 
\mathbb{G}_{\rho }(\langle \Sigma \rangle _{X})$. First, notice that one at
least of $\alpha x,\beta y$ is in $\Sigma .$ Otherwise, $1+\rho (\alpha
x)<0, $ as $x\in \Sigma $ and $\alpha x\in \langle x\rangle _{X}\backslash
\Sigma , $ and likewise $1+\rho (\beta y)<0.$ Summing,%
\begin{equation*}
2+\rho (\alpha x)+\rho (\beta y)<0.
\end{equation*}%
But $\alpha x+\beta y\in \mathbb{G}_{\rho }(X),$ so%
\begin{equation*}
0<1+\rho (\alpha x+\beta y)=1+\rho (\alpha x)+\rho (\beta y)<-1,
\end{equation*}%
a contradiction. We proceed by cases.

\noindent \textit{Case 1.} \textit{Both }$u:=\alpha x$ \textit{and} $%
v:=\beta y$\textit{\ are in} $\Sigma .$ Here 
\begin{equation*}
\alpha x+\beta y=u+v=u\circ _{\rho }[v/(1+\rho (u))]\in \Sigma ;
\end{equation*}%
indeed, by assumption $1+\rho (u))>0$ and $1+\rho (u+v)>0,$ so by linearity 
\begin{equation*}
1+\rho (v/(1+\rho (u)))=[1+\rho (u+v)]/(1+\rho (u))>0,
\end{equation*}%
and so $v/(1+\rho (u))\in \langle v\rangle _{\rho }\subseteq \Sigma .$

\noindent \textit{Case 2. One of} $u:=\alpha x,v=:\beta y$ \textit{is not in}
$\Sigma $\textit{\ (`off the half-line }$\langle x\rangle _{\rho }$\textit{\
or }$\langle y\rangle _{\rho }$\textit{').}

By commutativity of addition, without loss of generality (breiefly:
w.l.o.g.) $v\notin \Sigma .$ Then $-v\in \Sigma .$ As $\Sigma $ is a
subgroup, $(-v)^{-1}=v/(1-\rho (v))\in \Sigma $ and, setting%
\begin{equation*}
\delta :=(1-\rho (v))/[1+\rho (u)],
\end{equation*}%
\begin{equation*}
\alpha x+\beta y=u+v=u\circ _{\rho }\delta (-v)^{-1}=u+\delta v[1+\rho
(u)]/(1-\rho (v))\in \Sigma .
\end{equation*}%
Indeed, $\delta (-v)^{-1}=\delta v/(1-\rho (v))\in \langle v\rangle _{\rho
}\subseteq \Sigma ,$ since by assumption $1+\rho (u))>0$ and $1+\rho
(u+v)>0, $ so 
\begin{equation*}
1+\rho (\delta (-v)^{-1})=1+\rho \left( \frac{v}{1+\rho (u)}\right) =\frac{%
1+\rho (u+v)}{1+\rho (u)}>0.
\end{equation*}

Thus in all the possible cases $\alpha x+\beta y\in \Sigma $ for $x,y\in
\Sigma $ with $\alpha x+\beta y\in \mathbb{G}_{\rho }(\langle \Sigma \rangle
_{X}).$

Next we proceed by induction, with what has just been established as the
base step, to show that for all $n\geq 2,$ if $\alpha _{1}u_{1}+\alpha
_{2}u_{2}+...+\alpha _{n}u_{n}\in \mathbb{G}_{\rho }(\langle \Sigma \rangle
_{X}),$ for $u_{1},u_{2},...,u_{n}\in \Sigma $ and scalars $\alpha
_{1},\alpha _{2},...,\alpha _{n},$ then $\alpha _{1}u_{1}+\alpha
_{2}u_{2}+...+\alpha _{n}u_{n}\in \Sigma .$

Assuming the above for $n,$ we pass to the case of $u_{1},u_{2},...,u_{n+1}%
\in \Sigma $ and scalars $\alpha _{1},\alpha _{2},...,\alpha _{n+1}$ with $%
\alpha _{1}u_{1}+\alpha _{2}u_{2}+...+\alpha _{n+1}u_{n+1}\in \mathbb{G}%
_{\rho }(\langle \Sigma \rangle _{X}).$

Again as a preliminary, notice that, by permuting the subscripts as
necessary, w.l.o.g. $x:=\alpha _{1}u_{1}+...+\alpha _{n}u_{n}\in \mathbb{G}%
_{\rho }(\langle \Sigma \rangle _{X});$ otherwise, for $j=1,...,n+1$ 
\begin{equation*}
1+\rho \left( \sum\nolimits_{i\neq j}\alpha _{i}u_{i}\right) <0,
\end{equation*}%
and again as above, on summing, this leads to the contradiction%
\begin{equation*}
0<n[1+\rho (\alpha _{1}u_{1}+\alpha _{2}u_{2}+...+\alpha _{n+1}u_{n+1})]<-1.
\end{equation*}

So we suppose w.l.o.g. that $\alpha _{1}u_{1}+\alpha _{2}u_{2}+...+\alpha
_{n}u_{n}\in \mathbb{G}_{\rho }(\langle \Sigma \rangle _{X});$ by the
inductive hypothesis, $x:=\alpha _{1}u_{1}+\alpha _{2}u_{2}+...+\alpha
_{n}u_{n}\in \Sigma .$ Take $y:=u_{n+1}\in \Sigma $ and apply the base case $%
n=2$ to $x$ and $y.$ Then, since $w:=\alpha _{1}u_{1}+\alpha
_{2}u_{2}+...+\alpha _{n+1}u_{n+1}=x+\alpha _{n+1}y\in \mathbb{G}_{\rho
}(\langle \Sigma \rangle _{X})$, $w\in \Sigma .$ This completes the
induction, showing $\mathbb{G}_{\rho }(\langle \Sigma \rangle _{X})\subseteq
\Sigma .$ \hfill $\square $

\bigskip

In view of the role in quantifier weakening of countable subgroups dense in
themselves [BinO3,5], we note in passing that the proof above may be
relativized to the subfield of \textit{rational} scalars to give (with $%
\langle \cdot \rangle ^{_{\mathbb{Q}}}$ below the rational linear span):

\bigskip

\noindent \textbf{Theorem 1Q.} \textit{For }$\Sigma $\textit{\ a countable
subgroup of} $\mathbb{G}_{\rho }(X)$ \textit{with }$\langle u\rangle _{\rho
}^{\mathbb{Q}}\subseteq \Sigma $\textit{\ for each }$u\in \Sigma ,$ \textit{%
if} $\rho (\Sigma )\subseteq \mathbb{Q}:$%
\begin{equation*}
\Sigma =\mathbb{G}_{\rho }(\langle \Sigma \rangle ^{_{\mathbb{Q}}}).
\end{equation*}

\bigskip

\textbf{4. Abelian dichotomy and homomorphisms. }Our first result here,
Theorem 2, allows us to characterize in Theorems 4A and 4B homomorphisms
between Popa groups in vector spaces. We recall that%
\begin{equation*}
\eta _{1}(t):=1+t\qquad (t\in \mathbb{R}_{+}:=(0,\infty ))
\end{equation*}%
takes $\mathbb{G}_{1}(\mathbb{R)}\overset{\eta _{1}}{\rightarrow }(\mathbb{R}%
_{+},\times \mathbb{)}$, isomorphically. For the next result note that 
\begin{equation*}
\eta _{\rho }(x)=\eta _{1}(\rho (x))=1+\rho (x).
\end{equation*}%
In the case of $X=\mathbb{R}$, where $\rho (x)\equiv \rho x,$ this reduces to%
\begin{equation*}
1+\rho x.
\end{equation*}

We think of our first result here as expressing an \textit{abelian dichotomy}%
. Below $\circ _{I}$ refers to $\circ _{\rho }$ when $\rho =I,$ the identity
map on $\mathbb{R}$, as in the `circle operation' (above).

\bigskip

\noindent \textbf{Theorem 2.\newline
}\textit{A commutative subgroup }$\Sigma $\textit{\ of} $\mathbb{G}_{\rho
}(X)$ \textit{is either}\newline
\noindent (i)\textit{\ a subspace of the null space }$\mathcal{N}(\rho ),$ 
\textit{so a subgroup of }$(X,+),$\textit{\ or}\newline
\noindent (ii)\textit{\ for some }$u\in \Sigma $ \textit{a subgroup of }$%
\langle u\rangle _{\rho }$\textit{\ isomorphic under }$\rho $\textit{\ to a
subgroup of} $\mathbb{G}_{1}(\mathbb{R)}:$%
\begin{equation*}
\rho (x\circ _{\rho }y)=\rho (x)\circ _{I}\rho (y).
\end{equation*}%
\textit{Thus the image of }$\Sigma $\textit{\ under }$\eta _{\rho }$\textit{%
\ is a subgroup of }$(\mathbb{R}_{+},\times \mathbb{)}$.

\bigskip

\noindent \textbf{Proof.} Either $\rho (z)=0$ for each $z\in \Sigma ,$ in
which case $\Sigma $ is a subgroup of $(X,+),$ or else there is $z\in \Sigma
\backslash \{0\}$ with $\rho (z)\neq 0$ (since $\rho (0)=0$). In this case
take $u=u_{\rho }(z):=z/\rho (z)\neq 0.$ Then $\rho (u)=1$ so $u\in \Sigma ,$
and for all $x\in \Sigma $ by commutativity $x=\rho (u)x=\rho (x)u,$ i.e. $%
\Sigma $ is contained in the linear span $\langle u\rangle _{X}$ and so in $%
\langle u\rangle _{\rho }$. So the operation $\circ _{\rho }$ on $\Sigma $
takes the form%
\begin{equation*}
x\circ _{\rho }y=\rho (x)u+\rho (y)u+\rho (\rho (x)u)\rho (y)u.
\end{equation*}%
But $x\circ _{\rho }y=\rho (x\circ _{\rho }y)u,$ so as $u\neq 0$ the
asserted isomorphism follows from%
\begin{equation*}
\rho (x\circ _{\rho }y)u=[\rho (x)+\rho (y)+\rho (x)\rho (y)]u.
\end{equation*}%
In turn this implies%
\begin{equation*}
\eta _{\rho }(x\circ _{\rho }y)=1+\rho (x\circ _{\rho }y)=(1+\rho
(x))(1+\rho (y))=\eta _{\rho }(x)\eta _{\rho }(y),
\end{equation*}%
i.e. $\eta _{\rho }$ is a homomorphism into $(\mathbb{R}_{+},\times \mathbb{)%
}$.\hfill $\square $

\bigskip

Before we pass to a study of radial behaviours in \S 4, we recall the
following result [Ost2, Prop. A], [Chu1] (cf. [BinO5, Th. 3]) for the
context $\mathbb{G}_{\rho }(\mathbb{R})$ with $\rho (x)=\rho x.$ To
accommodate alternative forms of $(GFE),$ the matrix includes the
multiplicative group $(\mathbb{R}_{+},\times )$ as $\rho =\infty $ ; for a
derivation via a passage to the limit see [BinO5], but note that 
\begin{equation*}
\rho x+\rho y+\rho x\rho y=[\rho x\cdot \rho y](1+o(\rho ))\qquad (x,y\in 
\mathbb{R}_{+},\rho \rightarrow \infty ).
\end{equation*}

\bigskip

\noindent \textbf{Theorem BO.} \textit{Take }$\psi :\mathbb{G}_{\rho }(%
\mathbb{R}_{+})\rightarrow \mathbb{G}_{\sigma }(\mathbb{R})$\textit{\ a
homomorphism with }$\rho ,\sigma \in \lbrack 0,\infty ].$\textit{\ Then the
lifting }$\Psi :\mathbb{R}\rightarrow \mathbb{R}$\textit{\ \ of }$\psi $%
\textit{\ to }$\mathbb{R}$\ \textit{defined by the canonical isomorphisms }$%
\log ,\exp ,$\textit{\ }$\{\eta _{\rho }:\rho >0\}$ \textit{is bounded above
on }$\mathbb{G}_{\rho }$\textit{\ iff }$\Psi $\textit{\ is bounded above on }%
$\mathbb{R}$\textit{, in which case }$\Psi $ \textit{and }$\psi $\textit{\
are continuous. Then for some }$\kappa \in \mathbb{R}$ \textit{one has} $%
\psi (t)$ \textit{as below:}

\renewcommand{\arraystretch}{1.25}%
\begin{equation*}
\begin{tabular}{|l|l|l|l|}
\hline
Popa parameter & $\sigma =0$ & $\sigma \in (0,\infty )$ & $\sigma =\infty $
\\ \hline
$\rho =0$ & $\kappa t$ & $\eta _{\sigma }^{-1}(e^{\sigma \kappa t})$ & $%
e^{\kappa t}$ \\ \hline
$\rho \in (0,\infty )$ & $\log \eta _{\rho }(t)^{\kappa /\rho }$ & $\eta
_{\sigma }^{-1}(\eta _{\rho }(t)^{\sigma \kappa /\rho })$ & $\eta _{\rho
}(t)^{\kappa /\rho }$ \\ \hline
$\rho =\infty $ & $\log t^{\kappa }$ & $\eta _{\sigma }^{-1}(t^{\sigma
\kappa })$ & $t^{\kappa }$ \\ \hline
\end{tabular}%
\end{equation*}%
\medskip \newline
\renewcommand{\arraystretch}{1}

After linear transformation, all the cases reduce to some variant (mixing
additive or multiplicative structures) of the Cauchy functional equation.
(The parameters are devised to achieve continuity across cells, see [BinO5].)

\bigskip

We next show how this theorem is related to the current context of $(GFE)$.
As a preliminary we note a result of Chudziak in which $\circ _{\rho }$ is
applied to all of $X$, so in practice to Javor groups -- i.e. without
restriction to $\mathbb{G}_{\rho }(X).$ We thus think of this as a Javor
Homomorphism Theorem. We repeat Chudziak's proof, amending it to the range
context of $\mathbb{G}_{\sigma }(Y).$

\bigskip

\noindent \textbf{Theorem Ch} ([Chu2, Th. 1]). \textit{Let }$X,Y$ \textit{be
real topological vector spaces and } $K:X\rightarrow $ $\mathbb{G}_{\sigma
}(Y)$ \textit{a continuous function satisfying}%
\begin{equation*}
K(x\circ _{\rho }y)=K(x)\circ _{\sigma }K(y)\qquad (x,y\in X)
\end{equation*}%
\textit{with }$\rho \neq 0.$ \textit{Then for any }$u$\textit{\ with }$\rho
(u)=1$ \textit{there are constants }$\kappa =\kappa (u),\tau =\sigma (K(u)),$%
\textit{\ and continuous }$A_{u}:X\rightarrow $ $\mathbb{G}_{\sigma }(Y)$ 
\textit{satisfying}%
\begin{equation}
A_{u}(x+y)=A_{u}(x)\circ _{\sigma }A_{u}(y)\qquad (x,y\in X)  \tag{$A$}
\end{equation}%
\textit{(so with abelian range) such that}%
\begin{equation*}
K(x)=\left\{ 
\begin{array}{cc}
A_{u}(x)+[1+\sigma (A_{u}(x))][(1+\rho (x))^{\tau \kappa }-1]K(u)/\tau , & 
\tau \neq 0, \\ 
K(u)\log (1+\rho (x))/\log 2, & \tau =0.%
\end{array}%
\right.
\end{equation*}

\bigskip

\noindent \textbf{Proof. }Take any $u\in X$ with $\rho (u)=1$ and set%
\begin{equation*}
A_{u}(x):=K\left( x-\rho (x)u\right) ,\qquad \mu _{u}(t):=K((t-1)u).
\end{equation*}%
The former is continuous and satisfies $(A).$ To see this, take $%
v_{i}=x_{i}-\rho (x_{i})u;$ then $v_{1}+v_{2}=v_{1}\circ _{\rho }v_{2},$
since $\rho (v_{i})=\rho (x_{i})-\rho (x_{i})\rho (u)=0$ and $\circ _{\rho }$%
reduces to addition on the kernel of $\rho .$ Now, by linearity of $\rho ,$%
\begin{equation*}
v_{1}\circ _{\rho }v_{2}=v_{1}+v_{2}=x_{1}+x_{2}-\rho (x_{1}+x_{2})u.
\end{equation*}%
So%
\begin{eqnarray*}
A_{u}(x_{1}+x_{2}) &=&K\left( x_{1}+x_{2}-\rho (x_{1}+x_{2})u\right)
=K(v_{1}\circ _{\rho }v_{2}) \\
&=&K(v_{1})\circ _{\sigma }K(v_{2}) \\
&=&K\left( x_{1}-\rho (x_{1})u\right) \circ _{\sigma }K\left( x_{2}-\rho
(x_{2})u\right) \\
&=&A_{u}(x_{1})\circ _{\sigma }A_{u}(x_{2}).
\end{eqnarray*}%
Hence $A_{u}$ has image an abelian subgroup of $\mathbb{G}_{\sigma }(Y).$

The other mapping is an isomorphism between $(\mathbb{R}_{+},\times )$ and a
subgroup of $\mathbb{G}_{\sigma }(Y)$ with 
\begin{equation*}
\mu _{u}(st)=\mu _{u}(s)\circ _{\sigma }\mu _{u}(t).
\end{equation*}%
This last follows via $\rho (u)=1$ from the identity%
\begin{equation*}
(st-1)u=(s-1)u+[1+\rho ((s-1)u)](t-1)u.
\end{equation*}

Now the image subgroup under $\mu _{u}$, being abelian, is a subgroup of $%
\langle K(u)\rangle _{\sigma }$ by Theorem 2, so isomorphic to a subgroup of 
$\mathbb{G}_{\tau }(\mathbb{R})$ for $\tau :=\sigma (K(u))\in \mathbb{R}$.
Thus $\mu _{u}$ is an isomorphism from $(\mathbb{R}_{+},\times )=\mathbb{G}%
_{\infty }(\mathbb{R})$ to $\mathbb{G}_{\tau }(\mathbb{R}),$ for $\tau
=\sigma (K(u)),$ and by Theorem BO for some $\kappa =\kappa (u)$ 
\begin{equation*}
\mu _{u}(t)=\eta _{\sigma (K(u))}^{-1}(t^{\sigma (K(u))\kappa (u)})K(u).
\end{equation*}%
So, as $\rho ([x-\rho (x)u])=0,$ 
\begin{eqnarray*}
K(x) &=&K([x-\rho (x)u])\circ _{\rho }\rho (x)u)=A_{u}(x)\circ _{\sigma
}K(\rho (x)u) \\
&=&A_{u}(x)\circ _{\sigma }\mu _{u}(1+\rho (x)).
\end{eqnarray*}

For $\sigma (K(u))=0$ the above result should be amended to its limiting
value as $\tau \rightarrow 0,$ namely $K([x-\rho (x)u])+K(u)\log (1+\rho
(x))/\log 2$ (since $\kappa (u)=1/\log 2).$\hfill $\square $

\bigskip

\textbf{Remark.}\textit{\ }As the proof shows, in Theorem Ch. one fixes $u$\
with $\rho (u)=1,$\ obtaining constants $\kappa =\kappa (u),$ and $\tau
=\tau (u):=\sigma (K(u))$. The case $\tau =0$ is then best approached using
L'Hospital's rule so that, for $x=u,$ identity of both sides of the
representation of $K$ yields%
\begin{equation*}
1=\lim_{\tau \rightarrow 0}\frac{2^{\tau \kappa (u)}-1}{\tau }=\kappa
(u)\log 2.
\end{equation*}

\bigskip

\textbf{5. Radial behaviours. }Our next two results help establish in \S 6
Theorems 4A and 4B two not entirely dissimilar representations for the
circle groups, including the case $\rho \equiv 0,$ from which the form of $%
A_{u}$ above may be deduced in view of equation $(A)$ in Th. Ch. Our first
result concerns radial behaviour\textit{\ outside} $\mathcal{N}(\rho ).$

\bigskip

\noindent \textbf{Theorem 3A.}\newline
\textit{For real topological vector spaces }$X,Y$\textit{, if }$K:\mathbb{G}%
_{\rho }(X)\rightarrow \mathbb{G}_{\sigma }(Y)$ \textit{is continuous and
satisfies}%
\begin{equation}
K(x\circ _{\rho }y)=K(x)\circ _{\sigma }K(y)\qquad (x,y\in \mathbb{G}_{\rho
}(X)),  \tag{$K$}
\end{equation}%
\textit{then, for }$x$ \textit{with }$\rho (x)\neq 0$ \textit{and }$\sigma
(K(x))\neq 0,$\textit{\ there is }$\kappa =\kappa (x)\in \mathbb{R}$ $%
\backslash \{0\}$ \textit{with }%
\begin{equation*}
K(z)=\eta _{\sigma }^{-1}(\eta _{\rho }(z)^{\sigma (K(x))\kappa })\qquad
(z\in \langle x\rangle _{\rho }).
\end{equation*}%
\textit{Moreover, the index }$\gamma (x):=\sigma (K(x))\kappa (x)$\textit{\
is then continuous and extends to satisfy the equation}%
\begin{equation*}
\gamma (a\circ _{\rho }b)=\gamma (a)+\gamma (b)\qquad (a,b\in \mathbb{G}%
_{\rho }(X)).
\end{equation*}

\bigskip

\noindent \textbf{Proof.} For $x$ as above, take $u=u_{\rho }(x)\neq 0$ and $%
v=u_{\sigma }(K(x))\neq 0,$ both well-defined as $\rho (x)$ and $\sigma
(K(x))$ are non-zero (also $u$ $\in \langle x\rangle _{\rho }$ and $v$ $\in
\langle K(x)\rangle _{\sigma }$, as $\rho (u)=\sigma (v)=1$). The
restriction $K_{u}=K|\langle u\rangle _{\rho }$ yields a continuous
homomorphism into $\mathbb{G}_{\sigma }(Y).$ As $\langle u\rangle _{\rho }$
is an abelian group under $\circ _{\rho }$, its image under $K_{u}$ is an
abelian subgroup of $\mathbb{G}_{\sigma }(Y).$ So, as in Theorem 2, it is a 
\textit{non-trivial} subgroup of $\langle v\rangle _{\sigma }$. As noted, $%
\rho (u)=\sigma (v)=1,$ so we have the following \textit{isomorphisms}:%
\begin{eqnarray*}
&&\langle u\rangle _{\rho }\overset{\rho }{\rightarrow }\mathbb{G}_{1}(%
\mathbb{R)}\overset{\eta _{1}}{\rightarrow }(\mathbb{R}_{+},\times \mathbb{)}%
, \\
&&\langle v\rangle _{\sigma }\overset{\sigma }{\rightarrow }\mathbb{G}_{1}(%
\mathbb{R)}\overset{\eta _{1}}{\rightarrow }(\mathbb{R}_{+},\times \mathbb{)}
\end{eqnarray*}%
(writing $\rho ,\sigma =$ for $\rho |_{\langle u\rangle }$ and $\sigma
|_{\langle v\rangle })$ with $\langle .\rangle $ here short for $\langle
.\rangle _{\mathbb{R}}$), which combine to give%
\begin{equation*}
k(t):=\eta _{1}\sigma K_{u}\rho ^{-1}\eta _{1}^{-1}(t)=\eta _{\sigma
}K_{u}\eta _{\rho }^{-1}(t)
\end{equation*}%
as a \textit{non-trivial} homomorphism of $(\mathbb{R}_{+},\times \mathbb{)}$
into itself:%
\begin{equation*}
k(st)=k(s)k(t).
\end{equation*}%
Solving this Cauchy equation for a non-constant continuous $k$ yields 
\begin{equation*}
k(t)\equiv t^{\gamma }\qquad (t\in \mathbb{R}_{+}),
\end{equation*}%
for some $\gamma =\gamma (u)\in \mathbb{R}\backslash \{0\}$; so $k$ is
bijective. Write $\gamma =\gamma (u)=\sigma (K(u))\kappa (u)$, then, as
asserted (abbreviating the symbols),%
\begin{eqnarray*}
K_{u}(z) &=&\eta _{\sigma }^{-1}k\eta _{\rho }(z)=\eta _{\sigma }^{-1}(\eta
_{\rho }(z)^{\sigma \kappa }) \\
&=&\sigma ^{-1}(\eta _{1}^{-1}(1+\rho (z))^{\sigma \kappa }))\qquad (z\in
\langle u\rangle _{\rho }).
\end{eqnarray*}%
In particular, $K_{u}$ is injective. As $u\neq 0,$ $0\neq K(u)\in \langle
v\rangle _{\sigma },$ so $K(u)=sv$ for some $s\neq 0.$ Hence $\sigma
(K(u))=s\sigma (v)=s\neq 0.$ Since $\sigma (tK(u))=t\sigma (K(u)),$%
\begin{equation*}
K(z)=K_{u}(z)=[((1+\rho (z))^{\sigma (K(u))\kappa (u)}-1)/\sigma
(K(u))]K(u)\qquad (z\in \langle u\rangle _{\rho }).
\end{equation*}%
Here $\rho (z)=t$ for $z=tu,$ as $\rho (u)=1$ by choice. Taking $z=u$ gives%
\begin{equation*}
(2^{\sigma (K(u))\kappa (u)}-1)/\sigma (K(u))=1:\qquad \kappa (u)=\log
(1+\sigma (K(u))/[\sigma (K(u))\log 2],
\end{equation*}%
and so $\gamma (u):=\sigma (K(u))\kappa (u)$ is continuous and satisfies the
equation%
\begin{equation*}
\gamma (a\circ _{\rho }b)=\gamma (a)+\gamma (b)\qquad (a,b\in \mathbb{G}%
_{\rho }(X)).
\end{equation*}%
Indeed, write $\alpha =K(a),\beta =K(b);$ then as $K(a\circ _{\rho
}b)=\alpha \circ _{\sigma }\beta ,$ by linearity of $\sigma $%
\begin{eqnarray*}
\log (1+\sigma (K(a\circ _{\rho }b)) &=&\log (1+\sigma (\alpha +\beta
+\sigma (\alpha )\beta )) \\
&=&\log (1+\sigma (\alpha )+\sigma (\beta )+\sigma (a)\sigma (\beta )) \\
&=&\log (1+\sigma (\alpha ))+\log (1+\sigma (\beta )).
\end{eqnarray*}

For $\sigma (K(x))=0,$ the map $\langle v\rangle _{\sigma (v)}$ $\overset{%
\eta _{\sigma }}{\rightarrow }$ $(\mathbb{R}_{+},\times \mathbb{)}$ above
must be intepreted as exponential. A routine adjustment of the argument
yields%
\begin{equation*}
K(z)=K_{u}(z)=K(u)\log (1+\rho (z))/\log 2\qquad (z\in \langle u\rangle
_{\rho }),
\end{equation*}%
justifying hereafter a \textit{L'Hospital convention} (of taking limits $%
\sigma (K(u))\rightarrow 0$ in the `generic' formula).\hfill $\square $

\bigskip

We consider now the case $\rho (x)=0,$ which turns out as expected, despite
Theorem 2 being of no help here. This complement to Theorem 3A thus
describes radial behaviour\textit{\ inside} $\mathcal{N}(\rho )$. A more
detailed analysis, including the case $\rho (u)=1,$ along the lines followed
here, is to be found in [BinO7, Th. 3.1] and again in a Banach algebra
context in [BinO8].

\bigskip

\noindent \textbf{Theorem 3B.\newline
}\textit{Let }$X,Y$ \textit{be real topological vector spaces. If }$K:%
\mathbb{G}_{\rho }(X)\rightarrow \mathbb{G}_{\sigma }(Y)$\textit{\ is
continuous and satisfies }$(K),$\textit{\ then for any }$u\neq 0$\textit{\
with }$\rho (u)=0$\textit{\ }%
\begin{equation*}
K(\langle u\rangle _{\rho })\subseteq \langle K(u)\rangle _{\sigma }\sim 
\mathbb{G}_{\sigma (K(u))}\mathbb{(R)},
\end{equation*}%
\textit{and there is a function }$\lambda _{u}:(\mathbb{R},+)\rightarrow 
\mathbb{G}_{\sigma (K(u))}\mathbb{(R)}$ \textit{satisfying}%
\begin{equation*}
K(\xi u)=\lambda _{u}(\xi )K(u)\qquad (\xi \in \mathbb{R}).
\end{equation*}%
\textit{Moreover, if }$K(u)\neq 0,$\textit{\ then for some constant }$\kappa
=\kappa (u)$ 
\begin{equation*}
\lambda _{u}(t)=\left\{ 
\begin{array}{cc}
(e^{\sigma (K(u))\kappa (u)t}-1)/\sigma (K(u)), & \sigma (K(u))\neq 0, \\ 
t, & \sigma (K(u))=0,%
\end{array}%
\right.
\end{equation*}%
\textit{for }$t\in \mathbb{R}$, \textit{so that }$\lambda _{u}$ \textit{is
an isomorphism.}

\bigskip

\noindent \textbf{Proof}. As $\rho (u)=0,$ $\xi u+\xi u=\xi u\circ _{\rho
}\xi u.$ Notice that%
\begin{equation*}
K(2u)=K(u)+K(u)+\sigma (K(u))K(u)=(2+\sigma (K(u))K(u).
\end{equation*}%
By induction, 
\begin{equation*}
K(nu)=a_{n}(u)K(u)\in \langle K(u)\rangle _{Y},
\end{equation*}%
where $a_{1}=1$ and $a_{n}=a_{n}(u),$ for $n=1,2,...,$ solves%
\begin{equation*}
a_{n+1}=1+(1+\sigma (K(u))a_{n},
\end{equation*}%
since%
\begin{equation*}
K(u+nu)=K(u)+a_{n}K(u)(1+\sigma (K(u)).
\end{equation*}%
Suppose w.l.o.g. $\sigma (K(u))\neq 0,$ the case $\sigma (K(u))=0$ being
similar, but simpler (with $a_{n}=n).$ So%
\begin{equation*}
a_{n}=((1+\sigma (K(u))^{n}-1)/\sigma (K(u))\neq 0\qquad (n=1,2,...).
\end{equation*}%
Replacing $u$ by $u/n$ and then rearranging gives 
\begin{equation*}
K(u)=K(nu/n)=a_{n}(u/n)K(u/n):\qquad K(u/n)=a_{n}(u/n)^{-1}K(u)\in \langle
K(u)\rangle _{Y}.
\end{equation*}%
So%
\begin{eqnarray*}
K(mu/n) &=&a_{m}(u/n)K(u/n)=a_{m}(u/n)a_{n}(u/n)^{-1}K(u) \\
&=&\frac{((1+\sigma (K(u/n))^{m}-1)/\sigma (K(u/n)}{((1+\sigma
(K(u/n))^{n}-1)/\sigma (K(u/n)}K(u) \\
&=&\frac{((1+\sigma (K(u/n))^{m}-1)}{((1+\sigma (K(u/n))^{n}-1)}K(u)\in
\langle K(u)\rangle _{Y}.
\end{eqnarray*}

By continuity of $K$ (and of scalar multiplication), this implies that $%
K(\xi u)\in \langle K(u)\rangle _{Y}$ for any $\xi \in \mathbb{R}.$ So we
may uniquely define $\lambda (s)=\lambda _{u}(s)$ via%
\begin{equation*}
K(su)=\lambda _{u}(s)K(u).
\end{equation*}%
(In the case $\sigma (K(u))=0$ with $a_{n}=n,$ $K(mu/n)=(m/n)K(u),$ so that $%
K(su)=sK(u).)$ Then, as $\rho (u)=0,$%
\begin{eqnarray*}
\lambda (\xi +\eta )K(u) &=&K((\xi +\eta )u)=K(\xi u\circ _{\rho }\eta
u)=K(\xi u)+K(\eta u)+\sigma (K(\xi u))K(\eta u) \\
&=&K(\xi u)+K(\eta u)+\sigma (\lambda (\xi )K(u))\lambda (\eta )K(u) \\
&=&\lambda (\xi )K(u)+\lambda (\eta )K(u)+\lambda (\xi )\lambda (\eta
)\sigma (K(u))K(u) \\
&=&[\lambda (\xi )+\lambda (\eta )+\lambda (\xi )\lambda (\eta )\sigma
(K(u))]K(u).
\end{eqnarray*}%
So if $K(u)\neq 0$%
\begin{equation*}
\lambda _{u}(\xi +\eta )=\lambda _{u}(\xi )+\lambda _{u}(\eta )+\lambda
_{u}(\xi )\lambda _{u}(\eta )\sigma (K(u))=\lambda _{u}(\xi )\circ _{\sigma
(K(u))}\lambda _{u}(\eta ).
\end{equation*}%
Thus $\lambda _{u}:(\mathbb{R},+)\rightarrow \mathbb{G}_{\sigma (K(u))}(%
\mathbb{R}).$ By Theorem BO, with $\tau =\sigma (K(u))$ for some $\kappa
=\kappa (u)$%
\begin{equation*}
\lambda _{u}(t)=\left\{ 
\begin{array}{cc}
(e^{\tau \kappa (u)t}-1)/\tau , & \text{if }\sigma (K(u))\neq 0, \\ 
\kappa (u)t=t, & \text{if }\sigma (K(u))=0.%
\end{array}%
\right. \qquad \hfill \square
\end{equation*}

\bigskip

\noindent \textbf{Corollary 1.}\textit{\ In Theorem 3B, if }$\rho (u)=0$%
\textit{\ and }$K(u)\neq 0,$ \textit{then either}\newline
\noindent (i)\textit{\ }$\sigma (K(u))=0$ \textit{and }$\kappa (u)=1,$%
\textit{\ or\newline
}\noindent (ii) $\sigma (K(u))>0,$ $\kappa (u)=\log [1+\sigma (K(u))]/\sigma
(K(u))$ \textit{and the index }$\gamma (u):=\sigma (K(u))\kappa (u)$\textit{%
\ is additive on }$\mathcal{N}(\rho )$\textit{:}%
\begin{equation*}
\gamma (u+v)=\gamma (u)+\gamma (v)\qquad (u,v\in \mathcal{N}(\rho )).
\end{equation*}

\bigskip

\noindent \textbf{Proof. }As $\rho (u)=0,$ the notation in the proof above
is valid, so $\lambda _{u}(1)=1,$ as $0\neq K(u)=\lambda _{u}(1)K(u).$%
\textbf{\ }If \textit{\ }$\sigma (K(u))=0,$ then $\kappa (u)=1,$ by Theorem
3B. Otherwise,%
\begin{equation*}
(e^{\sigma (K(u))\kappa (u)}-1)/\sigma (K(u))=1:\qquad \kappa (u)=\log
(1+\sigma (K(u)))/\sigma (K(u)),
\end{equation*}%
and, as $\gamma (u)=\log (1+\sigma (K(u))),$ the concluding argument is as
in Theorem 3A (with $\circ _{\rho }=+$ on $\mathcal{N}(\rho )$). \hfill $%
\square $

\bigskip

\textbf{6. Homomorphism dichotomy. }The paired Theorems 4A and 4B below, our
main contribution, amalgamate the earlier radial results according to the
two forms identified by Theorem 2 that an \textit{abelian} Popa subgroup may
take (see below). Theorem 4A covers $\sigma \equiv 0$ as $\mathcal{N}(\sigma
)=\mathbb{G}_{\sigma }(Y)=Y,$ whereas $\rho \equiv 0$ may occur in the
context of either theorem. Relative to Theorem Ch., new here is Theorem 4B
exhibiting an additional source of regular variation.

\bigskip

We begin by noting that, since $\circ _{\rho }$ on $\mathcal{N}(\rho )$ is
addition, $\mathcal{N}(\rho )$ is an abelian subgroup of $\mathbb{G}_{\rho
}(X)$ and so 
\begin{equation*}
\Sigma :=K(\mathcal{N}(\rho )),
\end{equation*}%
as a homeomorph, is also an abelian subgroup of $\mathbb{G}_{\sigma }(Y).$
By Theorem 2 there are now two cases to consider, differing only in their
treatment of radial behaviour (in or out of $\mathcal{N}(\rho ))$. The
former is our First Popa Homomorphism Theorem\textit{\ }which follows.

\bigskip

\noindent \textbf{Theorem 4A.} \textit{Let }$X,Y\ $\textit{be real
topological vector spaces and }$K:\mathbb{G}_{\rho }(X)\rightarrow \mathbb{G}%
_{\sigma }(Y)$ \textit{a continuous function }$K$\textit{\ satisfying }$(K)$ 
\textit{with}%
\begin{equation*}
\mathit{\ }K(\mathcal{N}(\rho ))\subseteq \mathcal{N}(\sigma ).
\end{equation*}%
\textit{Then:}\newline
\noindent $K|\mathcal{N}(\rho )$\textit{\ is linear, and either}\newline
\noindent (i) $K$ \textit{is linear, or}

\noindent (ii)\textit{\ for any }$u$ \textit{with }$\rho (u)=1,\pi
_{u}(x):=x-\rho (x)u$ \textit{is the projection onto }$\mathcal{N}(\rho )$ 
\textit{parallel to }$u$ \textit{and }%
\begin{equation*}
K(x)=\left\{ 
\begin{array}{cc}
K(\pi _{u}(x))+K(u)[(1+\rho (x))^{\log (1+\tau )/\log 2}-1]/\tau , & \tau
\neq 0, \\ 
K(\pi _{u}(x))+K(u)\log (1+\rho (x))/\log 2, & \tau =0,%
\end{array}%
\right.
\end{equation*}%
\textit{for }$\tau =\sigma (K(u)).$\textit{\ In particular, }$x\mapsto K(\pi
_{u}(x))$\textit{\ is linear.}

\bigskip

\noindent \textbf{Proof.} If $\rho \equiv 0,$ then $K(X)=K(\mathcal{N}(\rho
))\subseteq \mathcal{N}(\sigma ).$ Here $\sigma (K(x))=0$ for all $x$ so,
since $\circ _{\sigma }=+$ on $\mathcal{N}(\sigma )$, $K$ is linear.

Otherwise, fix $u\in X$ with $\rho (u)=1.$ Then $x\mapsto \pi _{u}(x)=x-\rho
(x)u$ is a (linear) projection onto $\mathcal{N}(\rho )$ and, since $\rho
(x-\rho (x)u)=0,$%
\begin{equation*}
x=(x-\rho (x)u)\circ _{\rho }\rho (x)u.
\end{equation*}%
(So $\mathbb{G}_{\rho }(X)$ is generated by $\mathcal{N}(\rho )$ and any $%
u\notin \mathcal{N}(\rho ).$)

By assumption $\sigma (\pi _{u}(x))=0$ and as $K|\mathcal{N}(\rho )$ is
linear%
\begin{equation*}
K(x)=K(\pi _{u}(x))\circ _{\sigma }K(\rho (x)u)=K(\pi _{u}(x))+K(\rho (x)u).
\end{equation*}

If $\tau :=\sigma (K(u))\neq 0,$ then by Theorem 3A%
\begin{equation*}
K(\rho (x)u)=[(1+\rho (x))^{\log (1+\tau )/\log 2}-1]K(u)/\tau .
\end{equation*}

Now consider $u,v\in \mathbb{G}_{\rho }(X)$ with $\rho (u)=1=\rho (v).$ As $%
v-u\in \mathcal{N}(\rho ),$ also $\sigma (K(v-u))=0$ and also%
\begin{equation*}
v=(v-u)+u=(v-u)\circ _{\rho }u.
\end{equation*}%
Moreover, as $\sigma (K(v-u))=0,$%
\begin{equation*}
K(v)=K(v-u)\circ _{\sigma }K(u)=K(v-u)+K(u):\qquad K(v-u)=K(v)-K(u).
\end{equation*}%
So, by linearity of $\sigma $, 
\begin{equation*}
0=\sigma (K(v-u))=\sigma (K(v))-\sigma (K(u)):\qquad \sigma (K(v))=\sigma
(K(u))=\tau .
\end{equation*}%
Thus also%
\begin{equation*}
K(\rho (x)v)=[(1+\rho (x))^{\log (1+\tau )/\log 2}-1]K(v)/\tau .
\end{equation*}

If $\tau :=\sigma (K(u))=0,$ then as in Theorem 3A,%
\begin{equation*}
K(\rho (x)u)=K(u)\log (1+\rho (x))/\log 2,
\end{equation*}%
again justifying the L'Hospital convention in force (the formula follows
from the main case taking limits as $\tau \rightarrow 0$).\hfill $\square $

\bigskip

We pass to the remaining case, our Second Popa Homomorphism Theorem.

\bigskip

\noindent \textbf{Theorem 4B.} \textit{Let }$X,Y\ $\textit{be real
topological vector spaces and }$K:\mathbb{G}_{\rho }(X)\rightarrow \mathbb{G}%
_{\sigma }(Y)$ \textit{a continuous function }$K$\textit{\ satisfying }$(K)$ 
\textit{with }%
\begin{equation*}
K(\mathcal{N}(\rho ))=\langle K(w)\rangle _{\sigma }
\end{equation*}%
\textit{for some }$w$\textit{\ with }$\rho (w)=0$ \textit{and }$\sigma
(K(w))=1.$ \textit{Then}\newline
\noindent (i) $V_{0}:=\mathcal{N}(\rho )\cap K^{-1}(\mathcal{N}(\sigma ))$%
\textit{\ is a vector subspace and }$K_{0}=K|V_{0}=0$\textit{;\newline
\noindent }(ii) \textit{for any subspace }$V_{1}\ni $\textit{\ }$w$\textit{\
complementary to }$V_{0}$ \textit{in }$\mathcal{N}(\rho ),$\textit{\ and any 
}$u\in X$ \textit{with }$\rho (u)=1$\textit{, there is a linear map }$\kappa
_{w}:V_{1}\rightarrow \mathbb{R}$ \textit{such that for }$\tau =\sigma
(K(u)) $%
\begin{equation*}
K(x)=\left\{ 
\begin{array}{cc}
\lbrack e^{\kappa _{w}(\pi _{1}(x))}-1]K(w)+ &  \\ 
e^{\kappa _{w}(\pi _{1}(x))}[(1+\rho (x))^{\log (1+\tau )/\log
2}-1]K(u)/\tau , & \tau \neq 0, \\ 
\lbrack e^{\kappa _{w}(\pi _{1}(x))}-1]K(w)+ &  \\ 
e^{\kappa _{w}(\pi _{1}(x))}K(u)\log (1+\rho (x))/\log 2, & \tau =0,%
\end{array}%
\right.
\end{equation*}%
\textit{where }$\pi _{i}$\textit{\ denotes projection from }$X$\textit{\ onto%
} $V_{i}$\textit{, and }$\sigma (K(\pi _{1}(x)))\neq 0$ \textit{unless }$\pi
_{1}(x)=0$\textit{. }

\textit{(The final term in each case is excluded when there are no }$u\in X$ 
\textit{with }$\rho (u)=1$\textit{.)}

\bigskip

\noindent \textbf{Proof.} The assumption on $K\ $here is taken in an
initially more convenient form: $K(\mathcal{N}(\rho ))\subseteq \langle
w\rangle _{\sigma },$ \textit{for some} $w\in \Sigma =K(\mathcal{N}(\rho )),$
and of course w.l.o.g. $\sigma (w)\neq 0$, as otherwise this case is covered
by Theorem 4A.

To begin with $V_{0}:=\mathcal{N}(\rho )\cap K^{-1}(\mathcal{N}(\sigma ))$
is a subgroup of $\mathbb{G}_{\rho }(X),$ as $K$ is a homomorphism.
Similarly as in Theorem 4A, we work with a linear map, namely $%
K_{0}:=K|V_{0} $, as we claim $V_{0}$ to be a subspace of $\mathcal{N}(\rho
).$ (Then $V_{0}=\mathbb{G}_{0}(V_{0}).$)

The claim follows by linearity of $\sigma $ and Theorem 3B. Indeed, if $\rho
(x)=\rho (y)=0$ and $\sigma (K(x))=\sigma (K(y))=0,$ then $K(\alpha
x)=\lambda _{x}(\alpha )K(x)$ and $K(\beta y)=\lambda _{y}(\beta )K(y),$ and
since $\mathcal{N}(\rho )$ is a vector subspace on which $+$ agrees with $%
\circ _{\rho }:$%
\begin{eqnarray*}
K(\alpha x+\beta y) &=&K(\alpha x\circ _{\rho }\beta y) \\
&=&\lambda _{x}(\alpha )K(x)+\lambda _{y}(\beta )K(y)+\lambda _{x}(\alpha
)\lambda _{y}(\beta )\sigma (K(x))K(y) \\
&=&\lambda _{x}(\alpha )K(x)+\lambda _{y}(\beta )K(y),
\end{eqnarray*}
as $\sigma (K(x))=0.$ So%
\begin{equation*}
\sigma (K(\alpha x+\beta y))=\lambda _{x}(\alpha )\sigma (K(x))+\lambda
_{y}(\beta )\sigma (K(y))=0.
\end{equation*}%
Hence $V_{0}$ is a subspace of $\mathcal{N}(\rho )$ and $K_{0}:V_{0}%
\rightarrow \mathcal{N}(\sigma )$ is linear with $K_{0}(V_{0})\subseteq 
\mathcal{N}(\sigma ),$ as in Theorem 4A. Hence $K_{0}=0$; indeed, for $%
v_{0}\in V_{0},$ as $V_{0}\subseteq \mathcal{N}(\rho )$ there is $\lambda
_{0}$ with $K(v_{0})=\lambda _{0}w\in \mathcal{N}(\sigma )$ and so $0=\sigma
(\lambda _{0}w)=\lambda _{0}\sigma (w)$ and as $\sigma (w)\neq 0$ we have $%
\lambda _{0}=0.$ That is, $K_{0}=0.$

Since $K(\mathcal{N}(\rho ))\subseteq \mathcal{N}(\sigma )$ does not hold,
choose in $\mathcal{N}(\rho )$ a subspace $V_{1}$ complementary to $V_{0},$
and let $\pi _{i}:X\rightarrow V_{i}$ denote projection onto $V_{i}$. For $%
v\in \mathcal{N}(\rho )$ and $v_{i}=\pi _{i}(v)\in V_{i},$ as $K(v_{0})\in 
\mathcal{N}(\sigma ),$ 
\begin{equation}
K(v)=K(\pi _{0}(v)\circ _{\rho }\pi _{1}(v))=K(\pi _{0}(v))\circ _{\sigma
}K(\pi _{1}(v))=K_{0}(\pi _{0}(v))+K(\pi _{1}(v)).  \tag{$V0$}
\end{equation}%
Here $K_{0}\circ \pi _{0}$ is linear and $\sigma (K(v_{1}))\neq 0$ unless $%
v_{1}=0$. Recalling that $V_{1}$ is a subgroup of $\mathbb{G}_{\rho }(X),$
re-write the result of Theorem 3B as $K(v_{1})=\lambda _{w}(v_{1})w$ with $%
\lambda _{w}:V_{1}\rightarrow \mathbb{G}_{\sigma (w)}(\mathbb{R})$ and%
\begin{equation*}
\lambda _{w}(v_{1}+v_{1}^{\prime })=\lambda _{w}(v_{1})\circ _{\sigma
(w)}\lambda _{w}(v_{1}^{\prime })\qquad (v_{1},v_{1}^{\prime }\in V_{1}).
\end{equation*}%
With $w$ fixed, $\lambda _{w}$ is continuous (as $K$ is), with $1+\sigma
(w)\lambda _{w}(v_{1})>0$.

So, as in Theorem 4A, for $v\in V_{1}$ and some $\kappa =\kappa _{w}(v)$%
\begin{equation*}
K(tv)=\lambda _{w}(tv)w=\sigma (w)^{-1}[e^{\sigma (w)\kappa
_{w}(v)t}-1]w\qquad (t\in \mathbb{R}).
\end{equation*}

Taking $t=1$ gives%
\begin{equation*}
\sigma (w)\kappa _{w}(v)=\log [1+\sigma (w)\lambda _{w}(v)].
\end{equation*}%
As $\lambda _{w}$ is continuous, so is $\kappa _{w}:V_{1}\rightarrow \mathbb{%
R}$ (as $\sigma (w)\neq 0$). However, as in Theorem 2 but with $\sigma (w)$
fixed, $\kappa _{w}$ is additive and so by continuity linear on $V_{1}$. So,
as $t\kappa _{w}(v)=\kappa _{w}(tv),$ 
\begin{equation}
K(v)=\sigma (w)^{-1}[e^{\sigma (w)\kappa _{w}(v)}-1]w\qquad (v\in V_{1}). 
\tag{$V1$}
\end{equation}

For $x\in X$ take $v_{i}:=\pi _{i}(x)\in V_{i}$ and $v:=v_{0}+v_{1}.$ If $%
\rho $ is not identically zero, again fix $u\in X$ with $\rho (u)=1,$ and
then $x\mapsto \pi _{u}(x)=x-\rho (x)u$ is again (linear) projection onto $%
\mathcal{N}(\rho )$. If $\rho \equiv 0,$ set $u$ below to $0.$ Then, whether
or not $\rho \equiv 0,$ as $\rho (x-\rho (x)u)=0,$%
\begin{equation*}
x=v_{0}+v_{1}+\rho (x)u=v\circ _{\rho }\rho (x)u.
\end{equation*}%
So, as $\sigma (K(v_{0}))=0$ and $\rho (\rho (x)u)=\rho (x)\rho (u)=\rho
(x), $with $\tau =\sigma (K(u))\neq 0$%
\begin{equation*}
K(x)=K(v)\circ _{\sigma }K(\rho (x)u)=K(v)\circ _{\sigma }\eta _{\sigma
(K(u))}^{-1}(\eta _{\rho }(\rho (x)u)^{\kappa }),
\end{equation*}%
here with $\kappa =\log (1+\tau )/\log 2,$ which we expand as%
\begin{eqnarray*}
&&K(v_{0})+K(v_{1})+[1+\sigma (K(v_{0}+v_{1}))][(1+\rho (x))^{\log (1+\tau
)/\log 2}-1]K(u)/\tau \\
&=&K_{0}(\pi _{0}(v))+K(\pi _{1}(v))+[1+\sigma (K(v_{1}))][(1+\rho
(x))^{\log (1+\tau )/\log 2}-1]K(u)/\tau :
\end{eqnarray*}%
\begin{eqnarray*}
K(x) &=&K_{0}(\pi _{0}(x))+[e^{\sigma (w)\kappa _{w}(\pi
_{1}(x))}-1]w/\sigma (w) \\
&&+[1+\sigma (K(\pi _{1}(x)))][(1+\rho (x))^{\log (1+\sigma (K(u)))/\log
2}-1]K(u)/\sigma (K(u)).
\end{eqnarray*}%
Finally, $(V1)$ and linearity of $\sigma $ yields via $(V0)$ that%
\begin{equation*}
1+\sigma (K(\pi _{1}(x)))=e^{\sigma (w)\kappa _{w}(\pi _{1}(x))}.
\end{equation*}%
For $v_{1}\neq 0,$ $\sigma (K(v_{1}))\neq 0,$ as otherwise $v_{1}\in 
\mathcal{N}(\rho )\cap K^{-1}(\mathcal{N}(\sigma ))=V_{0},$ contradicting
complementarity of $V_{1}$.

Here $\sigma (w/\sigma (w))=1.$ Finally, as $w\in \Sigma =K(\mathcal{N}(\rho
)),$ we replace $w$ by $K(w)$ with $\rho (w)=0$ and $\sigma (K(w))=1.$ If $%
\tau =\sigma (K(u))=0,$ then, as in Theorem 4A, the final term is to be
interpreted by the L'Hospital convention (limiting value as $\tau =\sigma
(K(u))\rightarrow 0).$ We thus have:%
\begin{eqnarray*}
K(x) &=&K_{0}(\pi _{0}(x))+[e^{\kappa _{w}(\pi _{1}(x))}-1]K(w) \\
&&+e^{\kappa _{w}(\pi _{1}(x))}K(u)[(1+\rho (x))^{\log (1+\tau )/\log
2}-1]/\tau ,
\end{eqnarray*}%
where $\tau =\sigma (K(u)).$

If $\rho \equiv 0,$ then $u=0$ so that $K(u)=0,$ and the final term
vanishes.\hfill $\square $

\bigskip

\noindent \textbf{Remark. }We close with the final d\'{e}noument, which is
the connection between $(GFE)$ and Popa groups. Theorems 4A and 4B are used
in [BinO7] to characterize, for $X,Y$ real topological vector spaces, the
continuous solutions $K:\mathbb{G}_{\rho }(X)\rightarrow Y$ of $(GFE)$ as
homomorphisms between Popa groups $\mathbb{G}_{\rho }(X)$ and $\mathbb{G}%
_{\sigma }(Y)$ for some $\sigma .$ For an inkling of the context, notice that%
\textbf{\ }for\textit{\ }$K:\mathbb{G}_{\rho }(X)\rightarrow Y$\textit{\ }as
in Prop. 2.1,\ under the strong assumption that $K$\textit{\ }is injective,
a linear $\sigma :Y\rightarrow \mathbb{R}$ can readily be deduced yielding%
\begin{equation*}
K(u\circ _{\rho }v)=K(u)+g(u)K(v)=K(u)\circ _{\sigma }K(v),
\end{equation*}%
so that $K$ is a Popa homomorphism (cf. [Ost2, Th. 1]). We relax the strong
assumption in [BinO7].

\bigskip

\textbf{Acknowledgement. }We gratefully acknowledge the Referee's very
careful reading of our paper and the many wise and helpful suggestions for
improving clarity.

\bigskip

\textbf{References.}

\noindent \lbrack Bar] K. Baron, On the continuous solutions of the Go\l 
\k{a}b-Schinzel equation. \textsl{Aequationes Math.} \textbf{38} (1989), no.
2-3, 155--162.\newline
\noindent \lbrack Bin1] N. H. Bingham, Tauberian theorems and the central
limit theorem. \textsl{Ann. Prob.} \textbf{9} (1981), 221-231.\newline
\noindent \lbrack Bin2] N. H. Bingham, Scaling and regular variation. 
\textsl{Publ. Inst. Math. Beograd} \textbf{97} (111) (2015), 161-174.\newline
\noindent \lbrack Bin3] N. H. Bingham, Riesz means and Beurling moving
averages. \textsl{Risk and Stochastics} (Ragnar Norberg Memorial volume, ed.
P. M. Barrieu), Imperial College Press, 2019, Ch. 8, 159-172
(arXiv:1502.07494).\newline
\noindent \lbrack BinGT] N. H. Bingham, C. M. Goldie and J. L. Teugels, 
\textsl{Regular variation}, 2nd ed., Cambridge University Press, 1989 (1st
ed. 1987).\newline
\noindent \lbrack BinO1] N. H. Bingham and A. J. Ostaszewski, Homotopy and
the Kestelman-Borwein-Ditor theorem. \textsl{Canad. Math. Bull.} \textbf{54}
(2011), 12--20. \newline
\noindent \lbrack BinO2] N. H. Bingham and A. J. Ostaszewski, Beurling slow
and regular variation. \textsl{Trans. London Math. Soc. }\textbf{1} (2014)
29-56.\newline
\noindent \lbrack BinO3] N. H. Bingham and A. J. Ostaszewski, Cauchy's
functional equation and extensions: Goldie's equation and inequality, the Go%
\l \k{a}b-Schinzel equation and Beurling's equation. \textsl{Aequationes
Math.} \textbf{89} (2015), 1293--1310.\newline
\noindent \lbrack BinO4] N. H. Bingham and A. J. Ostaszewski, Beurling
moving averages and approximate homomorphisms. \textsl{Indag. Math. }\textbf{%
27} (2016), 601-633 (fuller version: arXiv1407.4093).\newline
\noindent \lbrack BinO5] {N. H. Bingham and A. J. Ostaszewski, }General
regular variation, Popa groups and quantifier weakening. \textsl{J. Math.
Anal. Appl.} \textbf{483} (2020) 123610, 31 pp. (arXiv1901.05996). \newline
\noindent \lbrack BinO6] {N. H. Bingham and A. J. Ostaszewski, Extremes and
regular variation. A lifetime of excursions through random walks and L\'{e}%
vy processes, 121--137, \textsl{Progr. Probab.}, 78, Birkh\"{a}%
user/Springer, Cham, 2021 (arXiv2001.05420).}\newline
\noindent \lbrack BinO7] {N. H. Bingham and A. J. Ostaszewski, }The Goldie
Equation: III. Homomorphisms from Functional Equations \textbf{(}initially
titled:\textbf{\ }Multivariate Popa groups and the Goldie Equation)
arXiv:1910.05817.\newline
\noindent \lbrack BinO8] {N. H. Bingham and A. J. Ostaszewski, }The Go\l 
\k{a}b-Schinzel and Goldie functional equations in Banach algebras,
arXiv:2105.07794.\newline
\noindent \lbrack BriD] N. Brillou\"{e}t and J. Dhombres, \'{E}quations
fonctionnelles et recherche de sous-groupes. \textsl{Aequationes Math.} 
\textbf{31} (1986), no. 2-3, 253--293.\newline
\noindent \lbrack Brz1] J. Brzd\k{e}k, Subgroups of the group Z$_{n}$ and a
generalization of the Go\l \k{a}b-Schinzel functional equation. \textsl{%
Aequationes Math. }\textbf{43} (1992), 59--71.\newline
\noindent \lbrack Brz2] J. Brzd\k{e}k, Bounded solutions of the Go\l \c{a}%
b-Schinzel equation. \textsl{Aequationes Math.} \textbf{59} (2000), no. 3,
248--254.\newline
\noindent \lbrack Chu1] J. Chudziak, Semigroup-valued solutions of the Go\l 
\k{a}b-Schinzel type functional equation. \textsl{Abh. Math. Sem. Univ.
Hamburg,} \textbf{76} (2006), 91-98.\newline
\noindent \lbrack Chu2] J. Chudziak, Semigroup-valued solutions of some
composite equations. \textsl{Aequationes Math.} \textbf{88} (2014), 183--198.%
\newline
\noindent \lbrack Chu3] J. Chudziak, Continuous on rays solutions of a Go\l 
\c{a}b-Schinzel type equation. \textsl{Bull. Aust. Math. Soc.} \textbf{91}
(2015), 273--277.\newline
\noindent \lbrack Coh] P. M. Cohn, \textsl{Algebra}, vol. 1, 2$^{\text{nd}}$
ed. Wiley, New York, 1982 (1st ed. 1974).\newline
\noindent \lbrack HewR] {E. Hewitt and K. A. Ross, \textsl{Abstract harmonic
analysis}. Vol. I, Grundl. math. Wiss. \textbf{115}, Springer 1963 [Vol. II,
Grundl. \textbf{152}, 1970].}\newline
\noindent \lbrack Jab1] E. Jab\l o\'{n}ska, Continuous on rays solutions of
an equation of the Go\l \c{a}b-Schinzel type. \textsl{J. Math. Anal. Appl.} 
\textbf{375} (2011), 223--229.\newline
\noindent \lbrack Jab2] E. Jab\l o\'{n}ska, Christensen measurability and
some functional equation. \textsl{Aequationes Math.} 81 (2011), 155--165.%
\newline
\noindent \lbrack Jac] N. Jacobson, \textsl{Lectures in Abstract Algebra}.
vol. I. Van Nostrand, New York, 1951.\newline
\noindent \lbrack Jav] P. Javor, On the general solution of the functional
equation $f(x+yf(x))=f(x)f(y)$. \textsl{Aequationes Math.} \textbf{1}
(1968), 235--238.\newline
\noindent \lbrack Ost1] A. J. Ostaszewski, Beurling regular variation, Bloom
dichotomy, and the Go\l \k{a}b-Schinzel functional equation. \textsl{%
Aequationes Math.} \textbf{89} (2015), 725-744. \newline
\noindent \lbrack Ost2] A. J. Ostaszewski, Homomorphisms from Functional
Equations: The Goldie Equation. \textsl{Aequationes Math. }\textbf{90}
(2016), 427-448 (arXiv: 1407.4089).\newline
\noindent \lbrack Pop] C. G. Popa, Sur l'\'{e}quation fonctionelle $%
f[x+yf(x)]=f(x)f(y).$ \textsl{Ann. Polon. Math.} \textbf{17} (1965), 193-198.%
\newline

\bigskip

\bigskip

\noindent Mathematics Department, Imperial College, London SW7 2AZ;
n.bingham@ic.ac.uk \newline
Mathematics Department, London School of Economics, Houghton Street, London
WC2A 2AE; A.J.Ostaszewski@lse.ac.uk\newpage

\end{document}